\documentclass[10pt,draft,twoside]{article}
\usepackage{amsmath,amsthm,amssymb}
\numberwithin{equation}{section}
\makeatletter
\def\@cite#1#2{{\upshape[{#1\if@tempswa , #2\fi }]}}
\let\@text=\text
\def\text#1{\@text{\upshape #1}}
\makeatother
\def\Z{\mathbb{Z}}
\def\Q{\mathbb{Q}}
\def\R{\mathbb{R}}
\def\C{\mathbb{C}}

\def\gU{\mathcal{H}} % generating function of U_d(n)
\def\gS{\mathcal{E}} % generating function of S(n;1^r)
\def\gP{\mathcal{S}} % generating function of S_k(n) ("Pamzv" no P desu)
\def\Re{\operatorname{Re}}

\def\Gal{\operatorname{Gal}}

\def\Li{Li}%{\operatorname{Li}}

\def\e#1{\varepsilon_{#1}}
\def\mzv#1{\zeta^\bullet_{#1}}

\def\finpamzv#1#2{S_{#1,#2}}
\def\kakko#1{\left(#1\right)}
\def\ckakko#1{\left\{#1\right\}}
\def\abs#1{\left|#1\right|}

\def\len#1{\ell(#1)}
\def\deq{:=}

\def\inprod#1#2{\left\langle#1,\,#2\right\rangle}

\newcommand{\wN}{\omega_{N}} % N-th root of unity
\newcommand{\wM}{\omega_{M}} % M-th root of unity
\def\tfrac#1#2{#1/#2} % slash fraction

\DeclareFontEncoding{OT2}{}{}
\DeclareFontSubstitution{OT2}{cmr}{m}{n}
\DeclareFontFamily{OT2}{cmr}{\hyphenchar\font45 }
\DeclareFontShape{OT2}{cmr}{m}{n}{%
   <5><6><7><8><9>gen*wncyr%
   <10><10.95><12><14.4><17.28><20.74><24.88>wncyr10}{}
\DeclareFontShape{OT2}{cmr}{b}{n}{%
   <5><6><7><8><9>gen*wncyb%
   <10><10.95><12><14.4><17.28><20.74><24.88>wncyb10}{}
\DeclareMathAlphabet{\mathcyr}{OT2}{cmr}{m}{n}
\DeclareMathAlphabet{\mathcyb}{OT2}{cmr}{b}{n}
\SetMathAlphabet{\mathcyr}{bold}{OT2}{cmr}{b}{n}
\def\LSh{L_{\scriptscriptstyle\mathcyr{Sh}\!}}
\def\LSt{L_*}

\newtheorem{thm}{Theorem}[section]
\newtheorem{lem}[thm]{Lemma}

\theoremstyle{definition}
\newtheorem{ex}[thm]{Example}
\newtheorem*{ackn}{Acknowledgement}

\theoremstyle{remark}
\newtheorem{rem}[thm]{Remark}

\newenvironment{keywords}{\smallskip\noindent{\bfseries Keywords:}}{}
\newenvironment{MSC}{\smallskip\noindent{\bfseries 2000 Mathematical Subject Classification:}}{}

\begin{document}
%*:;.;:*:;.;:*:;.;:*:;.;:*:;.;:*:;.;:*:;.;:*:;.;:*:;.;:*:;.;:*

%\title{\bfseries A partial alternating multiple zeta value}
%\title{\bfseries Partial alternating multiple zeta values\\
%arising from the non-commutative harmonic oscillator}
\title{\bfseries A variation of multiple $L$-values\\
arising from the spectral zeta function\\
of the non-commutative harmonic oscillator}
\author{Kazufumi KIMOTO%
\thanks{Partially supported by Grant-in-Aid for Young Scientists (B) No.\,20740021}
\ and Yoshinori YAMASAKI%
\thanks{Partially supported by Grant-in-Aid for JSPS Fellows No.\,19002485}}
\date{May 8, 2008}
%\date{December 12, 2007}
%\date{\today}
\maketitle

\begin{abstract}
A variation of multiple $L$-values,
%which we call the \emph{partial alternating multiple zeta value},
%and its generalizations.
%These sums
which arises from the description of
the special values of the spectral zeta function
of the non-commutative harmonic oscillator, is introduced.
In some special cases, 
we show that its generating function  
%generalized partial alternating multiple zeta values
can be written in terms of the gamma functions.
This result enables us to obtain
explicit evaluations of them. 
%the partial alternating multiple zeta values


\begin{MSC}
11M41, % Other Dirichlet series and zeta functions
05E05. % Symmetric functions
\end{MSC}
\end{abstract}

%=-=-=-=-=-=-=-=-=-=-=-=-=-=-=-=-=-=-=-=-=-=-=-=-=-=-=-=-=
\section{Introduction}
%=-=-=-=-=-=-=-=-=-=-=-=-=-=-=-=-=-=-=-=-=-=-=-=-=-=-=-=-=

The multiple zeta values
\begin{equation}
\mzv k(n_1,\dots,n_k)
\deq\sum_{1\le i_1<\dots< i_k}
\frac{1}{i_1^{n_1}i_2^{n_2}\dots i_k^{n_k}}
\end{equation} 
are natural extensions of the Riemann zeta value
$\zeta(n)=\sum_{i=1}^\infty i^{-n}$
introduced by Euler,
and have been of continuing interest to many mathematicians \cite{V}.
Recently, it has been shown by several authors that
they appear in various fields in mathematics
such as the knot invariant theory,
quantum group theory and mathematical physics
(see, e.g. \cite{KZ2001, Z1994}).
This fact implies the richness of the theory of the multiple zeta values
and encourages the recent studies of them.
%This fact encourages the recent studies of the multiple zeta values.
One of the main problems in studying multiple zeta values is 
to clarify the $\Q$-algebra structure
of the space spanned by them,
which is closely related to that of the category of mixed Tate motives.
In fact, for this purpose, a plenty of results concerning 
relations among them and  
exact calculations of them are investigated.
Furthermore, as a natural generalization,
Arakawa and Kaneko \cite{AK2004} introduce 
two kinds of multiple $L$-values
\begin{align}
\LSh(n_1,\dots,n_k;f_1,\dots,f_k)&\deq
\sum_{m_1>\dots>m_k>0}\frac{f_1(m_1-m_2)\dots f_{k-1}(m_{k-1}-m_k)f_k(m_k)}{m_1^{n_1}m_2^{n_2}\dots m_k^{n_k}},\label{eq:LSh}\\
\LSt(n_1,\dots,n_k;f_1,\dots,f_k)&\deq
\sum_{m_1>\dots>m_k>0}\frac{f_1(m_1)f_2(m_2)\dots f_k(m_k)}{m_1^{n_1}m_2^{n_2}\dots m_k^{n_k}},\label{eq:LSt}
\end{align}
where $f_1,\dots,f_k$ are $\C$-valued periodic functions on $\Z$
and also study their relations and exact evaluations.

%As in the case of the multiple zeta values,
%there are also many studies on
%various relations among the multiple $L$-values and
%exact evaluation of them.

In this paper,
we study the following variation $S^{(N,M)}_k(n_1,\ldots,n_k)$ 
($N,M\in\mathbb{N}$) of the multiple $L$-values;
\begin{equation}
\label{def:SNM}
S^{(N,M)}_k(n_1,\ldots,n_k)
\deq
\sum_{1\le i_1\le i_2\le\dots\le i_k}
\varepsilon^{(N)}_{i_1i_2\dots i_k}
\frac{\wM^{i_1+i_2+\dots+i_k}}{i_1^{n_1}i_2^{n_2}\dots i_k^{n_k}},
\end{equation}
where $\omega_M$ is a primitive $M$th root of unity and 
\begin{equation}
\varepsilon^{(N)}_{ij}\deq
\begin{cases}
0 & i=j\not\equiv0\pmod N\\
1 & \text{otherwise}
\end{cases}
=1-\delta_{ij}
%\kakko{1-\frac{1+\omega_N^{i}+\omega_N^{2i}+\cdots+\omega_N^{(N-1)i}}{N}}.
\kakko{1-\frac1N\sum_{r=0}^{N-1}\wN^{ri}},\qquad
\varepsilon^{(N)}_{i_1i_2\dots i_k}\deq\prod_{j=1}^{k-1}\varepsilon^{(N)}_{i_ji_{j+1}}.
\end{equation}
Here $\delta_{ij}$ is the Kronecker delta.
%Put $\varepsilon^{(N)}_{i_1i_2\dots i_k}\deq\prod_{j=1}^{k-1}\varepsilon^{(N)}_{i_ji_{j+1}}$.
%When $M=N$,
For simplicity, 
we sometimes write $S^{(N)}_k(n_1,\ldots,n_k)=S^{(N,N)}_k(n_1,\ldots,n_k)$, 
$S^{(N,M)}_k(n)=S^{(N,M)}_k(n,\ldots,n)$ and
$S^{(N)}_k(n)=S^{(N)}_k(n,\ldots,n)$.
%Furthermore when $k=1$,
%we understand that $\varepsilon^{(N)}_{i}\equiv 1$ so that  
We note that
$S^{(N,M)}_1(n)=\Li_{n}(\wM)$ where 
$\Li_n(z)\deq\sum^{\infty}_{i=1}\tfrac{z^{i}}{i^{n}}$ is 
the polylogarithm.

The aim of the paper is to establish generating function formulas
for the series $S_k^{(N,M)}(n)$, % when $M$ divides $N$,
and give an explicit evaluation of them 
in terms of Bernoulli numbers in the special case where $N=M=2$ and $n$ is even.
It is quite remarkable that the values $S_k^{(2)}(n)$ 
can be fully computable;
in fact, there are few examples of computable multiple $L$-values.
%though we have few examples among such kind of series. 
%(ordinary) multiple $L$-values.
% which are explicitly evaluated.
In this sense, 
$S_k^{(N)}(n)$ seems to be a nice variant of the ordinary multiple $L$-values.
%Though our intention is on the evaluation of the partial multiple $L$-values,
%it would be interesting to investigate the relations among them
%as well as the irrationality/transcendency of them.
%We leave these problems to the future study.

We will sometimes call $S^{(N,M)}_k(n_1,\ldots,n_k)$
as a \emph{partial multiple $L$-value}
because it is indeed a partial sum of the ``non-strict'' multiple $L$-value 
\begin{equation*}
%\label{for:N=1}
\sum_{1\le i_1\le i_2\le\dots\le
i_k}\frac{\wM^{i_1+i_2+\dots+i_k}}{i_1^{n_1}i_2^{n_2}\dots i_k^{n_k}}
=S^{(1,M)}_k(n_1,n_2,\dots,n_k).
\end{equation*}
In particular, $S^{(1)}_k(n_1,n_2,\dots,n_k)$ gives 
the non-strict multiple zeta value (see, e.g. \cite{Muneta2007}).
It is also worth remarking that $\e{ij}^{(N)}\to1-\delta_{ij}$ as $N\to\infty$ for fixed indices $i,j$,
so that we may regard the (strict) multiple $L$-values \eqref{eq:LSt}
as ``limiting case'' $S^{(\infty,M)}(n_1,n_2,\dots,n_k)$ of our series.
We notice that our partial multiple $L$-value $S^{(N,M)}_k(n_1,\ldots,n_k)$
is a special case of \emph{neither} the multiple $L$-values 
\eqref{eq:LSh} nor \eqref{eq:LSt}
since $\e{i_1 i_2\dots i_k}^{(N)}$ does depend on \emph{both}
the differences $i_j-i_{j-1}$ of adjacent indices
and the values of the indices $i_1,\dots,i_k$ themselves.
However, it is not difficult to see that 
$S^{(N,M)}_k(n_1,\ldots,n_k)$ can be expressed as a 
$\Q$-linear combination of \eqref{eq:LSh} (or \eqref{eq:LSt}).
Thus, for fixed $N$ and $M$,
it may be interesting to study the structure of the subalgebra 
spanned by all $S_k^{(N,M)}(n_1,n_2,\dots,n_k)$ 
in the $\Q$-algebra spanned by all multiple $L$-values
$S^{(1,M)}_k(n_1,\ldots,n_k)$.  
% \eqref{eq:LSh} (or \eqref{eq:LSt})
% with $N$-periodic functions $f_1, f_2, \dots$.
We leave these problems to the future study.

We now explain the spectral-theoretic origin of our series 
$S^{(N,M)}_k(n_1,\ldots,n_k)$.
A system of differential equations defined by the operator
\begin{align*}
Q\deq
\begin{pmatrix}\alpha & 0 \\ 0 & \beta\end{pmatrix}%
\left(-\frac12\frac{d^2}{dx^2}+\frac12x^2\right)%
+\begin{pmatrix}0 & -1 \\ 1 & 0\end{pmatrix}%
\left(x\frac{d}{dx}+\frac12\right)
\end{align*}
having two real parameters $\alpha, \beta$
is called the \emph{non-commutative harmonic oscillator}.
This system was first introduced and extensively studied
by Parmeggiani and Wakayama \cite{PW2001,PW2002} (see also \cite{P2008}).
It is shown that
when $\alpha,\beta>0$ and $\alpha\beta>1$,
$Q$ defines a positive, self-adjoint operator
on $L^2(\R)\otimes\C^2$ which has only a discrete spectrum
$
(0<)\lambda_1\le\lambda_2\le\dots\le\lambda_n\le\dots(\nearrow+\infty),
$
and the multiplicities of the eigenvalues are uniformly bounded.
In order to describe the total behavior of the spectrum,
Ichinose and Wakayama \cite{IW2005a} studied
the spectral zeta function
$\zeta_Q(s)\deq\sum_{n=1}^\infty\lambda_n^{-s}$
which is absolute convergent if $\Re(s)>1$.
This is analytically continued to the whole plane $\C$
and gives a single-valued meromorphic function
which has a simple pole at $s=1$
and `trivial' zeros at nonpositive even integers.
If $\alpha=\beta=1/\sqrt2$,
then $Q$ is unitarily equivalent to a couple of
the (ordinary) harmonic oscillators,
%\begin{align*}
%H=\kakko{-\frac12\frac{d^2}{dx^2}+\frac12x^2}
%\begin{pmatrix}
%1 & 0 \\ 0 & 1
%\end{pmatrix},
%\end{align*}
from which it follows that
$\zeta_Q(s)=2(2^s-1)\zeta(s)$.
%because
%$\Spec H=\ckakko{n+1/2\,;\,n=0,1,2,\dots}$
%with multiplicities exactly two.
Thus one can regard $\zeta_Q(s)$ as a deformation
of the Riemann zeta function $\zeta(s)$.

In describing the special values of the spectral zeta function $\zeta_Q(s)$,
the integrals
\[
 J_m(n)
=2^m\!\!\int_0^1\!\!\dots\!\int_0^1\!\!
\left(\!\frac{(1-x_1^4)(1-x_2^4\dotsb x_m^4)}
{(1-x_1^2\dotsb x_m^2)^2}\!\right)^{\!\!n\!\!}
\frac{dx_1\dotsb dx_m}{1-x_1^2\dotsb x_m^2}
\quad(m=2,3,4,\dots;\,n=0,1,2,\dots)
\]
and their generating functions
$g_m(x)=\sum_{n=0}^\infty\binom{-1/2}n J_m(n)x^n$
play a very important role.
%let $g_m(x)=\sum_{n=0}^\infty\binom{-1/2}n J_m(n)x^n$
%be a generating functions of $J_m(n)$. 
In fact, Ichinose and Wakayama \cite{IW2005b}
calculated the first two special values $\zeta_Q(2)$ and $\zeta_Q(3)$ 
in terms of $g_2(x)$ and $g_3(x)$, respectively.
The higher special values $\zeta_Q(m)$ ($m\ge4$)
are also expected to be
expressed by $g_m(x)$ and their generalizations 
%$m$-th Ap\'ery-like numbers $J_m(n)$
(see, e.g. \cite{O,KW2006a,K2007}).
%Thus it is interesting and important to study $J_m(n)$.
%
In the case where $m=2r$ is even, $J_{2r}(n)$ %the $m$-th Ap\'ery-like number
is explicitly given by
\begin{align*}
J_{2r}(n)=\sum_{p=0}^n(-1)^p\binom{-\frac12}p^{\!2}\binom np
\sum_{k=0}^{r-1}\zeta\!\kakko{2r-2k,\frac12}\finpamzv kp,
\end{align*}
where $\zeta(s,x)\deq\sum_{n=0}^\infty(n+x)^{-s}$
is the Hurwitz zeta function and
\begin{equation*}
\finpamzv kp=\sum_{1\le i_1\le i_2\le\dots\le i_k\le2p}\e{i_1i_2\dots i_k}^{(2)}
\frac{(-1)^{i_1+i_2+\dots+i_k}}{i_1^2i_2^2\dots i_k^2}.
\end{equation*}
Now it is immediate to see that
our series $S_k^{(N,M)}(n_1,\ldots,n_k)$ is a natural generalization 
of $S^{(2)}_k(2)=\lim_{p\to\infty}\finpamzv kp$
(we give the explicit formula of $S^{(2)}_k(2)$ in Example~\ref{example:Skeven}).

It is also worth remarking that
another kind of generating function
$w_2(t)=\sum_{n=0}^\infty J_2(n)t^n$ of $J_2(n)$
is regarded as a
period integral for the universal family of
the elliptic curves equipped with
a rational point of order $4$,
and satisfies a Picard-Fuchs differential equation
attached to this family of curves \cite{KW2006b}.

\subsection*{Conventions}
%---------------------------------------------------------

We recall several basic conventions on partitions and symmetric functions
(for further details, see \cite{M}).

A partition is a weakly decreasing sequence of nonnegative integers
which has finitely many nonzero entries.
For a partition $\lambda=(\lambda_1,\dots,\lambda_l)$ ($\lambda_l\ge1$),
the sum $\lambda_1+\dots+\lambda_l$
of entries in $\lambda$ is denoted by $\abs\lambda$
and the number $l$ of nonzero entries in $\lambda$ is denoted by $\len\lambda$.
We write $\lambda\vdash k$ to imply $\abs\lambda=k$,
and say $\lambda$ is a partition of $k$.
We denote by $\emptyset$ the (only) partition of $0$.
To indicate a multiple of the same numbers in $\lambda$,
we often write in an exponential form;
Let $m_i=m_i(\lambda)$ be the number of $i$'s in $\lambda$.
We call $m_i(\lambda)$ the multiplicity of $i$ in $\lambda$.
Then, we also write
$\lambda=(k^{m_k},\dots,2^{m_2},1^{m_1})$ or
$\lambda=1^{m_1}2^{m_2}\dots k^{m_k}$.
For instance, $\lambda=(4,2,2,1,1,1)$ is also written as
$\lambda=(4,2^2,1^3)=1^32^24^1$. 
When all the entries of $\lambda$ is even,
we call $\lambda$ an even partition.
For a given partition $\mu=(\mu_1,\dots,\mu_l)$ and a positive integer $q$,
we define $q\mu=(q\mu_1,\dots,q\mu_l)$.
We notice that
$\{\lambda\vdash 2k\,|\,\lambda\text{\,:\,even}\}
=\{2\mu\,|\,\mu\vdash k\}$.
If a given pair of two partitions $\lambda$ and $\mu$ satisfies that
$\lambda_i-\mu_i=0$ or $1$
for any index $i$,
then we say $\lambda/\mu$ is a vertical strip.
For instance, $(4,2,2,1,1,1)/(3,2,1,1)$ is a vertical strip.

Let $f(n)$ be a function on $\mathbb{N}$ and $a_n$ a sequence. 
Then, for a partition $\lambda$ and $q\in\mathbb{N}$,
we put $f(q\lambda):=\prod^{\ell(\lambda)}_{j=1}f(q\lambda_j)$
and $a_{q\lambda}:=\prod^{\ell(\lambda)}_{j=1}a_{q\lambda_j}$.
For instance, $(q\lambda)!=\prod^{\ell(\lambda)}_{j=1}(q\lambda_j)!$.  

Let $x_1,x_2,\dots$ be (infinitely many) variables.
For each positive integer $r$,
we respectively denote by $e_r=e_r(x_1,x_2,\dots)$ and
$h_r=h_r(x_1,x_2,\dots)$ the $r$-th elementary 
and $r$-th complete symmetric function defined by
\begin{align*}
e_r=\sum_{1\le i_1<i_2<\dots<i_r}x_{i_1}x_{i_2}\dots x_{i_r},\qquad 
h_r=\sum_{1\le i_1\le i_2\le\dots\le i_r}x_{i_1}x_{i_2}\dots x_{i_r}.
\end{align*}
We also put $e_0=h_0=1$ for convenience.
Moreover, for a partition $\lambda$, we put
$e_{\lambda}=\prod_{i\ge 1}e_{\lambda_i}$
and $h_{\lambda}=\prod_{i\ge 1}h_{\lambda_i}$.
The generating functions of $e_r$ and $h_r$ are given by
\begin{equation}
\label{eq:GFeh}
E(t)=\sum_{r=0}^\infty e_rt^r=\prod_{n=1}^\infty(1+x_nt),\qquad
H(t)=\sum_{r=0}^\infty h_rt^r=\prod_{n=1}^\infty(1-x_nt)^{-1}.
\end{equation}

%%%%%%%%%%%%%%%%%%%%%%%%%%%%%%%%%%%%%%%%%%%%%%%%%%%%%%%%%%%%%%%%%%%%%%%%%%%%%% 
%%%%%%%%%%%%%%%%%%%%%%%%%%%%%%%%%%%%%%%%%%%%%%%%%%%%%%%%%%%%%%%%%%%%%%%%%%%%%% 

\section{Generating functions}

% additional macros used in this section
\def\tS{\widetilde{S}^{(N,M)}}
\def\br{\boldsymbol{r}}
\def\bt{\boldsymbol{\tau}}
\def\be{\boldsymbol{\varepsilon}}
\def\then{\Rightarrow}
\def\excoef#1#2{C_{#1#2}}

In this section, 
we establish generating function formulas 
for the series $S^{(N,M)}_k(n)$.
To achieve this,
we first consider a decomposition
of the non-strict multiple sum $S^{(N,M)}_k(n)$
into the sum of several strict multiple sums.
Notice that each increasing sequence $1\le i_1\le i_2\le\dots\le i_k$ of
$k$ positive integers uniquely determines a sequence
$\br=(r_1,r_2,\dots,r_l)$, which we will refer to as the \emph{multiplicity} of the sequence
$(i_1,i_2,\dots,i_k)$, through the condition
\begin{align*}
\underbrace{i_1=\dots=i_{r_1}}_{r_1}
<\underbrace{i_{r_1+1}=\dots=i_{r_1+r_2}}_{r_2}
<i_{r_1+r_2+1}=\dots=i_{r_1+\dots+r_{l-1}}
<\underbrace{i_{r_1+\dots+r_{l-1}+1}=\dots=i_{r_1+\dots+r_l}}_{r_l}.
\end{align*}
Obviously, $\br$ is a permutation of a certain partition of $k$.
We denote by $\tS(n;\br)$ the partial sum of $S^{(N,M)}(n)$ whose running indices have multiplicity $\br$, i.e.
\begin{align*}
\tS(n;\br)=
\sum_{j_1<\dots<j_l}
\varepsilon^{(N)}_{\underbrace{\scriptstyle j_{1},\dots,j_{1}}_{r_1},%
\dots,\underbrace{\scriptstyle j_{l},\dots,j_{l}}_{r_l}}
\frac{\omega_M^{r_1j_{1}+\dots+r_lj_{l}}}
{j_{1}^{nr_1}\dots j_{l}^{nr_l}}
=\sum_{\substack{j_1<\dots<j_l\\ r_i>1\then N\mid j_i}}
\frac{\omega_M^{r_1j_{1}+\dots+r_lj_{l}}}
{j_{1}^{nr_1}\dots j_{l}^{nr_l}}.
\end{align*}
We also put
\begin{align*}
S^{(N,M)}(n;\emptyset)=1,\quad
S^{(N,M)}(n;\lambda)=\sum_{\br\in P(\lambda)}\tS(n;\br)\quad(\lambda\ne\emptyset),
\end{align*}
where $P(\lambda)$ denotes the set consisting of the permutations of a partition $\lambda$.
It is easy to see that
\begin{equation}
S^{(N,M)}_k(n)=\sum_{\lambda\vdash k}S^{(N,M)}(n;\lambda).
\end{equation}

\def\tR{\widetilde R^{(N,M)}_{p,q}(n;\br)}

To study the series $S^{(N,M)}(n;\lambda)$,
we here employ another function $R^{(N,M)}(n;\mu)$ 
defined by 
\begin{align*}
R^{(N,M)}(n;\emptyset):=1, \quad 
R^{(N,M)}(n;\mu)
\deq S^{(N,M)}(n;\mu_{>1})S^{(N,M)}(n;1^{m_1(\mu)})
\qquad (\mu\ne\emptyset).
\end{align*}
Here $\mu_{>1}$ denotes the partition defined by
$\mu_{>1}\deq2^{m_2(\mu)}3^{m_3(\mu)}\dots$.
%The following lemma is crucial. 
%Take a partition $\mu\vdash k$ such that $\mu_{>1}\ne\mu(\ne\emptyset)$ and
Fix a partition $\mu\vdash k$ and
put $q=m_1(\mu)$, $p=\len\mu-q$. We easily see that
\begin{align*}
R^{(N,M)}(n;\mu)=\sum_{\br\in P(\mu_{>1})}\tR,\qquad
\tR=\sum_{\substack{s_1<\dots<s_p\\ t_1<\dots<t_q\\ N\mid s_i}}
\frac{\wM^{r_1s_1+\dots+r_ps_p+t_1+\dots+t_q}}{s_1^{r_1n}\dots s_p^{r_pn}t_1^n\dots t_q^n}.
\end{align*}
In the sum $\tR$
%\begin{align*}
%\sum_{\substack{s_1<\dots<s_p\\ t_1<\dots<t_q\\ N\mid s_i}}
%\frac{\wM^{r_1s_1+\dots+r_ps_p+t_1+\dots+t_q}}{s_1^{r_1n}\dots s_p^{r_pn}t_1^n\dots t_q^n}
%\end{align*}
for each $\br\in P(\mu_{>1})$,
several of the running indices $t_1,\dots,t_q$ may coincide with certain $s_1,\dots,s_p$.
To describe the situation, we introduce the following map:
Put
\begin{align*}
I(p,q)=\left\{(\bt,\be)=(\tau_0,\tau_1,\dots,\tau_p,\e1,\dots,\e p)\,;\,
\tau_i\in\Z_{\ge0},\ \e i\in\{0,1\},\
\sum_{i=0}^p\tau_i+\sum_{i=1}^p\e i=q \right\}.
\end{align*}
For each element $(\br,(\bt,\be))\in P(\mu_{>1})\times I(p,q)$,
we associate a new sequence $\pi_\mu(\br,(\bt,\be))$ by
\begin{align*}
\pi_\mu(\br,(\bt,\be))=(1^{\tau_0},r_1+\e 1,1^{\tau_1},r_2+\e 2,\dots,r_p+\e p,1^{\tau_p}).
\end{align*}
Notice that there exists a partition $\lambda\vdash k$ such that
$\lambda/\mu_{>1}$ is a vertical strip and $\pi_\mu(\br,(\bt,\be))\in P(\lambda)$.
Namely, the correspondence $\pi_\mu$ defines a map
%\begin{align*}
$\pi_\mu:P(\mu_{>1})\times I(p,q)\to\coprod_{\substack{\lambda\vdash k\\ \lambda/\mu_{>1}:\text{vertical strip}}}P(\lambda).$
%\end{align*}
Thus it follows that
\begin{align*}
%R^{(N,M)}(n;\mu)&=
\sum_{\br\in P(\mu_{>1})}\tR
%\sum_{\substack{s_1<\dots<s_p\\ t_1<\dots<t_q\\ N\mid s_i}}
%\frac{\wM^{r_1s_1+\dots+r_ps_p+t_1+\dots+t_q}}{s_1^{r_1n}\dots s_p^{r_pn}t_1^n\dots t_q^n}
=\sum_{\substack{\lambda\vdash k\\ \lambda/\mu_{>1}:\text{vertical strip}}}
\sum_{\br\in P(\lambda)}\left|\pi_\mu^{-1}(\br)\right|\tS(n;\br).
\end{align*}
Since each $\left|\pi_\mu^{-1}(\br)\right|$ depends \emph{only} on $\lambda$,
we obtain
\begin{align*}
R^{(N,M)}(n;\mu)
&=\sum_{\substack{\lambda\vdash k\\ \lambda/\mu_{>1}:\text{vertical strip}}}
\left|\pi_\mu^{-1}(\lambda)\right|S^{(N,M)}(n;\lambda).
\end{align*}
Next, we calculate $\left|\pi_\mu^{-1}(\lambda)\right|$.
For each $a>2$, we assume that
$\lambda_{i_{a1}}=\dots=\lambda_{i_{a,d(a)}}=a$,
where $d(a)=m_a(\lambda)$.
Let us count the number of elements $(\br,(\bt,\be))$ in $I(p,q)$ such that
$\pi_\mu(\br,(\bt,\be))=\lambda$.
Notice that $\bt$ is uniquely determined by the assumption.
If $r_{i_{aj}}+\e{i_{aj}}=\lambda_{i_{aj}}=a$,
then it is possible that
$(r_{i_{aj}},\e{i_{aj}})=(a,0)\text{\ or\ }(a-1,1)$,
and there are exactly $\binom{m_a(\lambda)}{m_a(\lambda;\mu)}$ ways of the choice
of $i_{aj}$ such that $(r_{i_{aj}},\e{i_{aj}})=(a,0)$,
where $m_i(\lambda;\mu)=\left|\ckakko{j\,;\,\lambda_j=\mu_j=i}\right|$.
(Remark that $m_2(\lambda;\mu)=m_2(\lambda)$.)
Thus we have
$\left|\pi_\mu^{-1}(\lambda)\right|=\prod_{i>2}\binom{m_a(\lambda)}{m_a(\lambda;\mu)}$.
If $\mu$ is an \emph{even} partition and $\mu/\lambda_{>1}$ is a vertical strip,
then $m_i(\lambda;\mu)=m_i(\lambda)$ (if $i$ is even) or $0$ (if $i$ is odd) by definition,
and hence $\left|\pi_\mu^{-1}(\lambda)\right|=1$.
Consequently, we get the following lemma.
\begin{lem}
\label{lem:RN_by_SN}
For each $\mu\vdash k$, it holds that 
\begin{equation}
R^{(N,M)}(n;\mu)=\sum_{\substack{\lambda\vdash k\\ \lambda/\mu_{>1}
\text{:vertical strip}}}\prod_{i>2}\binom{m_i(\lambda)}{m_i(\lambda;\mu)} S^{(N,M)}(n;\lambda),
\end{equation}
where $m_i(\lambda;\mu)=\left|\ckakko{j\,;\,\lambda_j=\mu_j=i}\right|$.
In particular, if $\mu$ is even, then
\begin{equation}\label{eq:decomp_when_even}
R^{(N,M)}(n;\mu)=\sum_{\substack{\lambda\vdash k\\ \lambda/\mu_{>1}
\text{:vertical strip}}}S^{(N,M)}(n;\lambda).
\end{equation}
\qed
\end{lem}

\begin{lem}
\label{lem:N_unique_even}
For any $\lambda\vdash k$,
there uniquely exists $\mu\vdash k$ such that
$\mu_{>1}$ is even and $\lambda/\mu_{>1}$ is a vertical strip.
\end{lem}

\begin{proof}
It is immediate to see that
$\mu=1^{m_1(\lambda)+m_3(\lambda)+m_5(\lambda)+\dots}
2^{m_2(\lambda)+m_3(\lambda)}4^{m_4(\lambda)+m_5(\lambda)}
6^{m_6(\lambda)+m_7(\lambda)}\dots\vdash k$
is a unique partition
which satisfies all the desired conditions.  
%such that $\mu_{>1}$ is even, $\lambda/\mu_{>1}$ is a vertical strip.
%This proves the lemma.
%\begin{equation*}
%\mu=1^{m_1(\lambda)+m_3(\lambda)+m_5(\lambda)+\dots}
%2^{m_2(\lambda)+m_3(\lambda)}4^{m_4(\lambda)+m_5(\lambda)}
%6^{m_6(\lambda)+m_7(\lambda)}\dots.
%\end{equation*}
\end{proof}

By Lemmas \ref{lem:RN_by_SN} and \ref{lem:N_unique_even}, we readily obtain the
\begin{lem}
\label{lem:SN_by_sum}
Let $U^{(N,M)}_d(n)\deq \sum_{\mu\vdash d}S^{(N,M)}(n;2\mu)$.
Then it holds that 
\begin{equation}
\label{eq:NSkSU}
S^{(N,M)}_k(n)
=\sum_{\lambda\vdash k}S^{(N,M)}(n;\lambda)
=\sum_{\substack{\mu\vdash k\\ \mu_{>1}\text{:even}}}R^{(N,M)}(n;\mu)
=\sum_{0\le 2d\le k}S^{(N,M)}(n;1^{k-2d})U^{(N,M)}_d(n).
\end{equation}
\qed
\end{lem}
 
We next study the generating function of $S^{(N,M)}_k(n)$.
For this purpose, the following formula,
which is obtained by the canonical product expression of the gamma function,
is useful.

\begin{lem}
%[{\cite[Example~$11$, Chapter~$\mathrm{VI}$, p.\,$115$]{Bromwich1949}}]
\label{lem:infi-prod-formula}
For $a_{i},b_{i}\in\mathbb{C}$
satisfying $\sum^{l}_{i=1}a_i=\sum^{l}_{i=1}b_i$,
the equality
\begin{equation}\label{for:inf-prod}
\prod^{\infty}_{m=k}\prod_{j=1}^l\frac{m+a_j}{m+b_j}
=\prod_{j=1}^l\frac{\Gamma(k+b_j)}{\Gamma(k+a_j)}
%\prod^{\infty}_{m=k}\frac{(m+a_1)\cdots (m+a_l)}{(m+b_1)\cdots (m+b_l)}
%=\frac{\Gamma(k+b_1)\cdots \Gamma(k+b_l)}{\Gamma(k+a_1)\cdots \Gamma(k+a_l)}
\end{equation}
holds for any integer $k$. \qed
\end{lem}

\begin{lem}
\label{lem:GF4U}
The generating function of $U^{(N,M)}_d(n)$
%{\upshape(}which depends only on $N${\upshape)}
is given by
\begin{equation}
\label{for:gen_UN}
\gU^{(N,M)}(n;x)\deq\sum_{d=0}^\infty U^{(N,M)}_d(n)x^{2nd}
%=\prod^{2n-1}_{j=0}\Gamma\!\kakko{1-\frac{x}{N}\omega_{2n}^j}.
=\prod_{k=1}^M\prod_{j=0}^{2n-1}
\frac{\Gamma\!\kakko{\frac1M\!\kakko{k-\frac1N\omega_{2n}^j\omega_{Mn}^{kN}x}}}
{\Gamma\!\kakko{\frac kM}}.
\end{equation}
\end{lem}

\begin{proof}
We notice that
\begin{align*}
U^{(N,M)}_d(n)=\sum_{\mu\vdash d}S^{(N)}(n;2\mu)
=h_d\Bigl(\frac{\wM^{2N}}{N^{2n}},\frac{\wM^{4N}}{(2N)^{2n}},\frac{\wM^{6N}}{(3N)^{2n}},\ldots\Bigl)
\end{align*}
since the complete symmetric function $h_d$ is
the sum of \emph{all} monomials of degree $d$.
Therefore, by specializing $x_m=\tfrac{\wM^{2mN}}{(Nm)^{2n}}$ and $t=x^{2n}$
in the generating function $H(t)$ in \eqref{eq:GFeh}, we obtain 
\begin{align*}
\gU^{(N,M)}(n;x)
&=\prod^{\infty}_{m=1}\kakko{1-\frac{\wM^{2mN}}{(Nm)^{2n}}x^{2n}}^{\!\!-1}
=\prod_{m=0}^\infty\prod_{k=1}^M\ckakko{1-\kakko{\frac{\omega_{Mn}^{kN}x}{N(Mm+k)}}^{\!2n}}^{\!\!-1}\\
&=\prod_{m=0}^\infty\prod_{k=1}^M\prod_{j=0}^{2n-1}
\kakko{1-\omega_{2n}^j\frac{\omega_{Mn}^{kN}x}{N(Mm+k)}}^{\!\!-1}
=\prod_{m=0}^\infty\prod_{k=1}^M\prod_{j=0}^{2n-1}
\frac{m+\frac kM-\frac{\omega_{2n}^j\omega_{Mn}^{kN}x}{MN}}{m+\frac kM}.
%&=\prod^{\infty}_{m=1}\prod^{2n-1}_{j=0}\kakko{1-\omega_{2n}^j\frac{x}{Nm}}^{\!\!-1}
%=\prod^{\infty}_{m=1}\prod^{2n-1}_{j=0}\frac{m}{m-\frac{x}{N}\omega_{2n}^j}.
\end{align*}
Applying Lemma \ref{lem:infi-prod-formula} to the equation above,
we have \eqref{for:gen_UN}.
%by using the equation \eqref{for:inf-prod}.
\end{proof}

\begin{lem}
\label{lem:GF4S}
The generating function of $S^{(N,M)}(n;1^r)$
%{\upshape(}which depends only on $M${\upshape)}
is given by
\begin{equation}
\label{for:gen_S1}
\gS^{(M)}(n;x)\deq\sum_{r=0}^\infty S^{(N,M)}(n;1^r)x^{nr}
=\prod_{k=1}^M\prod_{j=0}^{n-1}
\frac{\Gamma\!\kakko{\frac{k}{M}}}
{\Gamma\!\kakko{\frac1M\!\kakko{k-\omega^{2j-1}_{2n}\omega_{Mn}^{k}x}}}.
\end{equation}
\end{lem}

\begin{proof}
We notice that
\begin{align*}
S^{(N)}(n;1^r)
=e_r\Bigl(\frac{\wM}{1^n},\frac{\wM^2}{2^n},\frac{\wM^3}{3^n},\ldots\Bigr).
\end{align*}
Hence, if we specialize $x_m=\tfrac{\wM^m}{m^n}$ and 
set $t=x^n$ in the generating function $E(t)$ in \eqref{eq:GFeh}, 
then we obtain the lemma by a similar calculation as in the case of $\gU^{(N,M)}(n;x)$.
%\begin{align*}
%\gS^{(M)}(n;x)
%&=\prod_{m=1}^{\infty}\kakko{1+\frac{\wM^m}{m^n}x^n}
%=\prod_{m=0}^{\infty}\prod^{M}_{k=1}\kakko{1+\frac{\wM^{k}}{(Mm+k)^n}x^{n}}\\
%&=\prod_{m=0}^{\infty}\prod^{M}_{k=1}\prod_{j=0}^{n-1}\kakko{1-\omega_n^j\frac{\omega_{2Mn}^{2k-M}}{Mm+k}x}
%=\prod_{m=0}^{\infty}\prod^{M}_{k=1}\prod_{j=0}^{n-1}\frac{m+\frac kM-\frac xM\omega_{2n}^{2j-1}\omega_{Mn}^{k}}{m+\frac kM}
%\kakko{1-\omega_n^j\frac{\omega_{2Nn}^{N+2k}}{Nm+k}x}
%&=\prod_{m=1}^{\infty}\prod^{N-1}_{k=0}\kakko{1+\frac{\omega_N^{(Nm-k)}}{(Nm-k)^n}x^{n}}
%=\prod_{m=1}^{\infty}\prod^{N-1}_{k=0}\kakko{1+\frac{\omega_N^{-k}}{(Nm-k)^n}x^{n}}\\
%&=\prod_{m=1}^{\infty}\prod^{N-1}_{k=0}\prod^{n-1}_{j=0}\kakko{1-\omega^{j}_{n}\frac{\omega_{2Nn}^{N-2k}}{Nm-k}x}
%=\prod_{m=1}^{\infty}\frac{\prod^{N-1}_{k=0}\prod^{n-1}_{j=0}\kakko{m-\frac{k}{N}-\frac{x}{N}\omega^{j}_{n}\omega_{2Nn}^{N-2k}}}{\prod^{N-1}_{k=0}\kakko{m-\frac{k}{N}}^n}.
%\end{align*}
%Hence we obtain the desired formula by using \eqref{for:inf-prod} again.
%Hence, 
%using the equation \eqref{for:inf-prod} and rewriting $k$ as $N-k$,
%we obtain the desired formula. 
\end{proof}

Now, we obtain the following 
\begin{thm}\label{thm:result1}
The generating function of $S^{(N,M)}_k(n)$ is given by
\begin{equation}\label{for:gen_SNk_general}
\gP^{(N,M)}(n;x)\deq\sum_{k=0}^\infty S^{(N)}_k(n)x^{nk}
=\prod_{k=1}^M
\frac{\prod_{j=0}^{2n-1}\Gamma\!\kakko{\frac1M\!\kakko{k-\frac1N{\omega_{2n}^j\omega_{Mn}^{kN}x}}}}
{\Gamma\!\kakko{\frac kM}^{\!n}\prod_{j=0}^{n-1}\Gamma\!\kakko{\frac1M\!\kakko{k-{\omega_{2n}^{2j-1}\omega_{Mn}^{k}x}}}}.
%=\prod_{k=1}^M
%\frac{\prod_{j=0}^{2n-1}\Gamma\!\kakko{\frac kM-\frac{\omega_{2n}^j\omega_{Mn}^{kN}x}{MN}}}
%{\Gamma\!\kakko{\frac kM}^n\prod_{j=0}^{n-1}\Gamma\!\kakko{\frac kM-\frac{\omega_{n}^j\omega_{2Mn}^{2k-M}x}{M}}}.
%=\frac{\prod^{M}_{k=1}\Gamma\!\kakko{\frac{k}{M}}^n\prod^{2n-1}_{j=0}\Gamma\!\kakko{1-\frac{x}{N}\omega_{2n}^j}}
%{\prod^{M}_{k=1}\prod^{n-1}_{j=0}\Gamma\!\kakko{\frac{k}{M}-\frac{x}{M}\omega^{j}_{n}\omega_{2Mn}^{2k-M}}}.
\end{equation}
\end{thm}
\begin{proof}
From the equation \eqref{eq:NSkSU}, it is clear that 
$\gP^{(N,M)}(n;x)=\gU^{(N,M)}(n;x)\gS^{(M)}(n;x)$.
Hence one immediately obtains the formula \eqref{for:gen_SNk_general}
from \eqref{for:gen_UN} and \eqref{for:gen_S1}.
\end{proof}

If $M\mid N$, then, using the Gauss-Legendre formula of the gamma function,
we have the following reduced formulas:
\begin{align}
\gU^{(N,M)}(n;x)
&=\prod^{2n-1}_{j=0}\Gamma\!\kakko{1-\frac{\omega_{2n}^jx}{N}},\\
\gP^{(N,M)}(n;x)
&=\prod_{k=1}^M
\frac{\Gamma\!\kakko{\frac{k}{M}}^n\prod^{2n-1}_{j=0}\Gamma\!\kakko{1-\frac{\omega_{2n}^jx}{N}}}
{\prod^{n-1}_{j=0}\Gamma\!\kakko{\frac1M\!\kakko{k-\omega^{2j-1}_{2n}\omega_{Mn}^{k}x}}}.\label{for:gen_SNk}
\end{align}
Notice that $\gS^{(M)}(n;x)$ depends \emph{only} on $M$.

%In the next section, we give an explicit formula of
%the partial alternating multiple zeta value $S_k^{(2)}(n)$
%by using the result obtained above.

\section{Partial alternating multiple zeta values}

In this section,
we concentrate on the special case where $N=M=2$.
From the definition,
the sums $S_k(n)\deq S_k^{(2,2)}(n)$ in this case
may be called \emph{partial alternating multiple zeta values}. 
%$S_k(n)\deq S_k^{(2,2)}(n)
%=\sum_{1\le i_1\le i_2\le\dots\le i_k}\e{i_1i_2\dots i_k}\frac{(-1)^{i_1+i_2+\dots+i_k}}{i_1^ni_2^n\dots i_k^n}$
%which, from the definition,
%We give an explicit formula of
%the \emph{partial alternating multiple zeta value} 
%By using the result on generating functions obtained above.
From \eqref{for:gen_SNk}, we have 
\[
\gP(n;x):=\gP^{(2)}(n;x)
=\frac{\Gamma\!\kakko{\frac{1}{2}}^n\prod^{2n-1}_{j=0}\Gamma\!\kakko{1-\frac{x}{2}\omega_{2n}^{j}}}{\prod^{n-1}_{j=0}\Gamma\!\kakko{\frac{1}{2}-\frac{x}{2}\omega_{n}^j}\Gamma\!\kakko{1-\frac{x}{2}\omega_{n}^j\omega_{2n}}}
=\frac{\Gamma\!\kakko{\frac{1}{2}}^n\prod^{n-1}_{j=0}\Gamma\!\kakko{1-\frac{x}{2}\omega_{n}^{j}}}{\prod^{n-1}_{j=0}\Gamma\!\kakko{\frac{1}{2}-\frac{x}{2}\omega_{n}^{j}}}.
\]
Furthermore, using the duplication formula
$\Gamma\kakko{2a}\Gamma\!\kakko{\tfrac{1}{2}}=2^{2a-1}\Gamma(a)\Gamma\!\kakko{\tfrac{1}{2}+a}$
with $a=-\tfrac{x\omega_{n}^{j}}{2}$
and the equation $\sum^{n-1}_{j=0}\omega_{n}^{j}=\delta_{n,1}$, 
we see
\begin{equation}
\label{for:gen_S2k} 
\gP(n;x)
=\frac{\Gamma\!\kakko{\frac{1}{2}}^n\prod^{n-1}_{j=0}\Gamma\!\kakko{1-\frac{x}{2}\omega_{n}^{j}}}{\prod^{n-1}_{j=0}\Gamma(-x\omega_{n}^{j})\Gamma\!\kakko{\frac{1}{2}}2^{x\omega_n^j+1}\Gamma\!\kakko{-\frac{x}{2}\omega_{n}^{j}}^{-1}}
=2^{-x\delta_{n,1}}\prod^{n-1}_{j=0}\frac{\Gamma\!\kakko{1-\frac{x}{2}\omega_{n}^{j}}^2}{\Gamma\!\kakko{1-x\omega_{n}^{j}}}.
\end{equation} 

For $m\ge 0$, define the sequence $\{A^{\bullet}(m)\}_{m\ge 0}$ 
by $A^{\bullet}(0):=1$, $A^{\bullet}(1):=0$ and 
\[
A^{\bullet}(m):= 
\sum_{a=1}^{m-1}\mzv a(\underbrace{1,\dots,1}_{a-1},m-a+1) \qquad (m\ge 2).
\]
Namely, $A^{\bullet}(m)$ ($m\ge2$) denotes 
the sum of multiple zeta values of weight $m$ and height $1$.
It is known that $A^{\bullet}(m)$ can be expressed as a polynomial in
$\zeta(2),\zeta(3),\ldots$ and $\zeta(m)$ with rational coefficients 
(see \cite{OZ2001}).
For example, we have
$A^{\bullet}(3)=\zeta(3)+\mzv2(1,2)=2\zeta(3)$
since $\mzv2(1,2)=\zeta(3)$, which is due to Euler.
Further, we put
\begin{align*}
A^{\bullet}_{n}(m)&:=
\sum_{\substack{m_1,\ldots,m_{n}\ge 0\\m_1+\cdots+m_{n}=m}}
A^{\bullet}(m_1)\cdots A^{\bullet}(m_{n}),\qquad
 Z_n(k):=
\sum_{\substack{\mu\vdash k\\\mu_{\ell(\mu)}>\delta_{n,1}}}
\frac{\nu(n\mu)}{z_{\mu}}\zeta(n\mu),
\end{align*}
where $\nu(x):=2^{1-x}-1$ and $z_{\mu}:=\prod_{i\ge1}i^{m_i(\mu)}m_{i}(\mu)!$. 
Note that $A^{\bullet}_{1}(m)=A^{\bullet}(m)$. 
Then, we get the following expressions of the values 
$S_k(n)=S^{(2)}_k(n)$.

\begin{thm}
$(\mathrm{i})$\ If $n=1$, then it holds that
\begin{equation}
\label{for:n=1}
S_k(1)
=\sum^{k}_{m=0}\frac{(-\log{2})^{k-m}}{(k-m)!2^m}A^{\bullet}(m)
=\sum^{k}_{m=0}\frac{(-\log{2})^{k-m}}{(k-m)!}Z_{1}(m)
\in\mathbb{Q}[\,\log{2},\zeta(2),\zeta(3),\ldots,\zeta(k)].  
\end{equation} 
$(\mathrm{ii})$\ If $n\ge 2$, then it holds that
\begin{equation}
\label{for:nge2}
S_k(n)
=\frac{1}{2^{nk}}A^{\bullet}_{n}(nk)=Z_{n}(k)
\in\mathbb{Q}[\zeta(n),\zeta(2n),\ldots,\zeta(kn)]. 
\end{equation}
\end{thm} 
\begin{proof}
From the generating function \eqref{for:gen_S2k},
it is sufficient to show that 
\begin{equation}
\label{for:key}
 \prod^{n-1}_{j=0}\frac{\Gamma\!\kakko{1-\frac{x}{2}\omega_{n}^{j}}^2}{\Gamma\!\kakko{1-x\omega_{n}^{j}}}
=\sum^{\infty}_{m=0}A^{\bullet}_{n}(nm)\kakko{\frac{x}{2}}^{nm}
=\sum^{\infty}_{m=0}Z_n(m)x^{nm}.
\end{equation}
To prove this, we recall the identity
(see \cite{Aomoto1990,Drinfeld1991})
\begin{align}
\label{eq:ADZ}
\frac{\Gamma(1-X)\Gamma(1-Y)}{\Gamma(1-X-Y)}
=1-\sum_{a,b=1}^\infty
\mzv{a}(\underbrace{1,\dots,1}_{a-1},b+1)X^aY^b
=\exp\kakko{\sum_{m=2}^{\infty}\frac{X^m+Y^m-(X+Y)^m}{m}\zeta(m)}.
\end{align}
Putting $X=Y=\tfrac{x\omega_{n}^j}{2}$ 
and writing $a+b=m$ in the middle term in \eqref{eq:ADZ}, 
we have
\[
 \frac{\Gamma\!\kakko{1-\frac{x}{2}\omega_{n}^{j}}^2}{\Gamma\!\kakko{1-x\omega_{n}^{j}}}
=\sum^{\infty}_{m=0}A^{\bullet}(m)\kakko{\frac{x\omega_{n}^{j}}{2}}^m
=\exp\kakko{\sum_{m=2}^\infty\frac{\nu(m)}{m}\zeta(m)(\omega_n^jx)^m}.
\]
Then, taking the product $\prod^{n-1}_{j=0}$ of this equation, 
one sees that 
\begin{equation}
\label{for:middle-key}
 \prod^{n-1}_{j=0}\frac{\Gamma\!\kakko{1-\frac{x}{2}\omega_{n}^{j}}^2}{\Gamma\!\kakko{1-x\omega_{n}^{j}}}
=\sum^{\infty}_{m=0}A^{\bullet}_n(m)\kakko{\frac{x\omega_{n}^{j}}{2}}^m
=\exp\kakko{\sum_{\substack{m=1\\ nm\ge 2}}^{\infty}\frac{\nu(nm)}{m}\zeta(nm)x^{nm}}
\end{equation}
because $\sum_{j=0}^{n-1}\omega_n^{jm}=n$ if $n\mid m$ and $0$ otherwise.
Here, the rightmost-hand side of \eqref{for:middle-key}
can be written as
\begin{align*}
 \prod_{\substack{m=1\\ nm\ge 2}}^{\infty}\exp\kakko{\frac{\nu(nm)}{m}\zeta(nm)x^{nm}}
&=\prod_{\substack{m=1\\ nm\ge
 2}}^{\infty}\sum^{\infty}_{l_m=0}\frac{1}{l_m!}\kakko{\frac{\nu(nm)}{m}\zeta(nm)x^{nm}}^{l_m}\\
&=
\begin{cases}
 \displaystyle{\sum_{l_2,l_3,\ldots=0}^{\infty}\frac{\nu(2)^{l_2}\nu(3)^{l_3}\cdots}{(2^{l_2}3^{l_3}\cdots)(l_2!l_3!\cdots)}\kakko{\zeta(2)^{l_1}\zeta(3)^{l_3}\cdots}x^{2l_2+3l_3+\cdots}} & \textrm{($n=1$)}\\
 \displaystyle{\sum_{l_1,l_2,\ldots=0}^{\infty}\frac{\nu(n)^{l_1}\nu(2n)^{l_2}\cdots}{(1^{l_1}2^{l_2}\cdots)(l_1!l_2!\cdots)}\kakko{\zeta(n)^{l_1}\zeta(2n)^{l_2}\cdots}x^{l_1+2l_2+\cdots}} & \textrm{($n\ge 2$)}
\end{cases}\\
&=\sum^{\infty}_{m=0}\Biggl\{
\sum_{\substack{\mu\vdash m\\ \mu_{\ell(\mu)}>\delta_{n,1}}}
\frac{\nu(n\mu)}{z_{\mu}}\zeta(n\mu)\Biggr\}x^{nm}
=\sum^{\infty}_{m=0}Z_n(m)x^{nm}.
\end{align*}
Note that, from the second equality in \eqref{for:middle-key}, 
this shows that $A^{\bullet}_n(m)=0$ if $n\nmid m$. 
Therefore, one can actually obtain the equations \eqref{for:key}.
This completes the proof of the theorem.
\end{proof}

\begin{ex}
We have
\begin{align*}
S_1(1)=-\log2,\quad
S_2(1)=\frac{(\log2)^2}2-\frac{\zeta(2)}4,\quad
S_3(1)=-\frac{(\log2)^3}6+\frac{\log2}4\zeta(2)-\frac14\zeta(3),
\end{align*}
and
\begin{align*}
S_1(3)=-\frac34\zeta(3),\quad
S_2(3)=-\frac{31}{64}\zeta(6)+\frac9{32}\zeta(3)^2,\quad
S_3(3)=-\frac{255}{768}\zeta(9)+\frac{93}{128}\zeta(6)\zeta(3)-\frac{27}{384}\zeta(3)^3.
\end{align*}
%\begin{align*}
%S_1(n)=(2^{1-n}-1)\zeta(n),\quad
%S_2(n)=\frac{2^{1-2n}-1}2\zeta(2n)+\frac{(2^{1-n}-1)^2}2\zeta(n)^2,\\
%S_3(n)=\frac{2^{1-3n}-1}3\zeta(3n)+\frac{(2^{1-2n}-1)(2^{1-n}-1)}2\zeta(2n)\zeta(n)+\frac{(2^{1-n}-1)^3}6\zeta(n)^3.
%\end{align*}
\end{ex}

If one further assumes that $n$ is even,
then one can obtain the following various expressions.

\begin{thm}
It holds that
\begin{align}
S_k(2n)
&=(-\pi^2)^{nk}
\sum_{\substack{m_1,\ldots,m_{n}\ge 0\\ m_1+\cdots+m_{n}=nk}}
\omega_n^{m_1+2m_2+\dots+nm_{n}}
\frac{B_{2m_1}}{(2m_1)!}\cdots\frac{B_{2m_{n}}}{(2m_{n})!}
\label{eq:S2n}\\
&=(-\pi^2)^{nk}
\sum_{\substack{\lambda\vdash nk\\\len\lambda\le n}}
\inprod{p_n\circ h_k}{m_\lambda}\frac{B_{2\lambda}}{(2\lambda)!}
\label{eq:S2n2}\\
&=(-\pi^2)^{nk}\sum_{\mu\vdash k}
\frac{\widetilde{\nu}(2n\mu)}{z_\mu}\frac{B_{2n\mu}}{(2n\mu)!},
\label{eq:S2n3}
\end{align}
where $\widetilde{\nu}(x):=2^{x-1}-1$, 
$p_n$ is the $n$-th power-sum symmetric function,
$m_\lambda$ the monomial symmetric function for $\lambda$,
$\circ$ the plethysm,
and $\inprod\cdot\cdot$ the standard scalar product
in the ring of symmetric functions
defined by $\inprod{h_{\lambda}}{m_{\mu}}=\delta_{\lambda\mu}$
with $\delta_{\lambda\mu}$ being the Kronecker delta
$($see \cite{M} for detail\/$)$.
\end{thm}
\begin{proof}
If we apply the reflection formula
for the gamma function in \eqref{for:gen_S2k},
then we have
\begin{equation}\label{for:S2nk}
\begin{split}
\gP(2n;x)
&=\prod_{j=1}^{n}
\frac{\Gamma\!\kakko{1-\frac{\omega_{2n}^jx}{2}}^{\!2}\Gamma\!\kakko{1+\frac{\omega_{2n}^jx}{2}}^{\!2}}
{\Gamma\!\kakko{1-\omega_{2n}^jx}\Gamma\!\kakko{1+\omega_{2n}^jx}}
=\prod_{j=1}^{n}
%\kakko{
\frac{\pi x\omega_{2n}^j}{2}\cot\frac{\omega_{2n}^j\pi x}2
%}
=\prod_{j=1}^{n}
%\ckakko{
\sum_{m=0}^\infty
\frac{(-\omega_n^j)^mB_{2m}\pi^{2m}}{(2m)!}x^{2m},
%},
\end{split}
\end{equation}
from which we immediately obtain \eqref{eq:S2n}.
Next, it readily follows from \eqref{eq:S2n} that
\begin{align*}
S_k^{(2n)}=(-\pi^2)^{nk}
\sum_{\substack{\lambda\vdash nk\\\len\lambda\le n}}
m_\lambda(1,\omega_n,\dots,\omega_n^{n-1},0,\dots)
\frac{B_{2\lambda}}{(2\lambda)!}.
\end{align*}
Thus we should calculate
$m_\lambda(1,\omega_n,\dots,\omega_n^{n-1},0,\dots)$.
Let us recall the expansion formula (see, e.g. \cite{M})
\begin{align}
\label{eq:Cauchy}
\prod_{i,j\ge1}\frac1{1-x_iy_j}
=\sum_{\lambda}h_\lambda(x)m_\lambda(y).
\end{align}
If we set $y_j=\omega_n^{j-1}$ for $j=1,2,\dots,n$ and $y_j=0$ for $j>n$
in \eqref{eq:Cauchy}, then we have
\begin{align}\label{eq:plethysm}
\sum_{\len\lambda\le n}h_\lambda(x)m_\lambda(1,\omega_n,\dots,\omega_n^{n-1},0,\dots)
=\prod_{i\ge1}\frac1{1-x_i^n}
=\sum_{k=0}^\infty h_k(x_1^n,x_2^n,\ldots)
=\sum_{k=0}^\infty p_n\circ h_k.
\end{align}
By taking the terms of homogeneous degree $nk$ in \eqref{eq:plethysm},
we have
\begin{align}\label{eq:plethsym}
\sum_{\substack{\lambda\vdash nk\\ \len\lambda\le n}}
h_\lambda m_\lambda(1,\omega_n,\dots,\omega_n^{n-1},0,\dots)
=p_n\circ h_k
\end{align}
for each $k$.
Hence we get
$m_\lambda(1,\omega_n,\dots,\omega_n^{n-1},0,\dots)=\inprod{p_n\circ h_k}{m_\lambda}$,
which readily implies \eqref{eq:S2n2}.
The equation \eqref{eq:S2n3} follows immediately from \eqref{for:nge2}
together with the classical result $\zeta(2m)=(-1)^{m-1}2^{2m-1}B_{2m}\pi^{2m}/(2m)!$
due to Euler. %\eqref{eq:EulerZeta}.
This completes the proof.
\end{proof}

%We note that $m_\lambda(1,\omega_n,\dots,\omega_n^{n-1},0,\dots)\in\mathbb{Z}$
%since $\{h_\lambda\}_\lambda$ forms a $\Z$-basis
%of the ring of symmetric functions.

\begin{ex}
\label{example:Skeven}
From the equation \eqref{eq:S2n}, we have 
\begin{align*}
 S_{k}(2)
=\frac{(-1)^kB_{2k}}{(2k)!}{\pi}^{2k}=-\frac{\zeta(2k)}{2^{2k-1}},\qquad 
 S_{k}(4)
=\biggl\{\sum^{2k}_{m=0}(-1)^m\frac{B_{2m}B_{4k-2m}}{(2m)!(4k-2m)!}\biggr\}\pi^{4k}.
\end{align*}
%In particular, one obtains Theorem~\ref{thm:main}.
See \cite{Y2007}
for a similar discussion on the multiple Dirichlet $L$-values.
\end{ex}

\begin{rem}
It is remarkable that $S_k(2)=S_k^{(2)}(2)$ can be reduced as above.
% This simplicity of the result for $S_k(2)$ is remarkable.
We recall that $S_k^{(2)}(2)$ is closely related to 
the special value $\zeta_Q(2)$ of the spectral zeta function.
%$S_{k}^{(2)}(2)=-\zeta(2k)/2^{2k-1}$.
Can one explain the simplicity (or ``exact solvability'') of $S_k^{(2)}(2)$ by,
for instance, the existence of the Picard-Fuchs differential equation for $w_2(t)$?
%the modular property of $w_2(t)$?
%We expect that $\finpamzv kp$ and its limit $S_k^{(2)}(2)$ reflect
%a certain property of the special value $\zeta_Q(k)$ and/or attached generating functions $g_k(x)$, $w_k(t)$.
%For instance,

%the limit value $S_k^{(2)}(2)$ is explicitly given by
%\begin{equation}\label{eq:mainresult}
%S_k^{(2)}(2)=\frac{(-1)^{k}B_{2k}\pi^{2k}}{(2k)!},
%\end{equation}
%where $B_m$ denotes the $m$-th Bernoulli number (see Section 3).
%because of the following classical result due to Euler (see, e.g. \cite{V}):
%\begin{equation}\label{eq:EulerZeta}
%\zeta(2m)=\frac{(-1)^{m-1}2^{2m-1}B_{2m}\pi^{2m}}{(2m)!}.
%\end{equation}
\end{rem}

\begin{rem}
Let us give an example of the partial alternating {\it double} zeta value
with distinct indices:
\begin{align*}%\label{eq:Euler-like}
S^{(2)}_2(1,2k)&=(k+1)S^{(2)}_1(2k+1)+2(1-2^{-2k})\zeta(2k)\log2%
-\sum_{p=1}^{k-1}S^{(2)}_1(2p+1)\zeta(2k-2p),\\
%S_2(1,2k)=(k+1)\LSh(2k+1;\chi)+(2-2^{1-2k})\zeta(2k)\log2%
%+\sum_{p=1}^k\LSh(2p+1;\chi)\zeta(2k-2p),
S^{(2)}_2(2k,1)&=-kS^{(2)}_1(2k+1)-\zeta(2k)\log2%
+\sum_{p=1}^{k-1}S^{(2)}_1(2p+1)\zeta(2k-2p).
%S_2(1,2k)=(k+1)\LSh(2k+1;\chi)+(2-2^{1-2k})\zeta(2k)\log2%
%+\sum_{p=1}^k\LSh(2p+1;\chi)\zeta(2k-2p),
\end{align*}
%where $\chi(n)=(-1)^n$
Notice that $S^{(2)}_1(n)=(2^{1-n}-1)\zeta(n)$ for $n\ge2$.
This is regarded as an analogue of Euler's formula
$\mzv2(1,2k)=k\zeta(2k+1)-\frac12\sum_{p=2}^{2k-1}\zeta(p)\zeta(2k-p+1)$.
%To prove \eqref{eq:Euler-like},
%we have to evaluate the double $L$-value $\LSh(2k,1;\chi,1)$
%since $S_2(1,2k)=\LSh(2k,1;\chi,1)+\zeta(2k+1)/2^{2k+1}$.
%Actually, it is calculated by using the shuffle relations for multiple $L$-values (see \cite{DLV2008}).
See also \cite{BZB2008} for related calculations.
\end{rem}

\begin{ackn}
The authors would like to thank Professor Masato Wakayama
for valuable comments.
\end{ackn}

%/%/%/%/%/%/%/%/%/%/%/%/%/%/%/%/%/%/%/%/%/%/%/%/%/%/%/%/%/%/%/%/%/%/

%/%/%/%/%/%/%/%/%/%/%/%/%/%/%/%/%/%/%/%/%/%/%/%/%/%/%/%/%/%/%/%/%/%/

\bigskip

\noindent
\textsc{Kazufumi KIMOTO}\\
Department of Mathematical Sciences,
University of the Ryukyus\\
Senbaru, Nishihara, Okinawa 903-0231, Japan\\
\texttt{kimoto@math.u-ryukyu.ac.jp}

\medskip

\noindent
\textsc{Yoshinori YAMASAKI}\\
Faculty of Mathematics,
Kyushu University\\
Hakozaki, Fukuoka 812-8581, Japan\\
\texttt{yamasaki@math.kyushu-u.ac.jp}

%*:;.;:*:;.;:*:;.;:*:;.;:*:;.;:*:;.;:*:;.;:*:;.;:*:;.;:*:;.;:*
\end{document}